# CLASSIFICATION OF RECURRING UNEMPLOYED WORKERS AND UNEMPLOYMENT EXITS

Marie Cottrell et Patrice Gaubert[1]

*SAMOS-MATISSE et SET-MATISSE*
*Université Paris 1*
*90, rue de Tolbiac, 75634 PARIS Cedex 13*
*E-mail : cottrell, gaubert@univ-paris1.fr*

**Abstract** : This study focuses on recurring unemployment, that is people with two or more spells of unemployment during the period of observation (July 1993 – August 1996). First, a classification is obtained which is then used to examine the specific role of occasional jobs during a spell of unemployment and, in this context, the influence of the received unemployment benefits on the duration of this spell. This paper is a continuation of previous analyses of unemployment in France, based on long-term data from the unemployed register held by ANPE (National Employment Bureau). The present analysis conducted using additional information about unemployment benefits received by the unemployed from UNEDIC (Unemployment Benefits Office).

## 1. INTRODUCTION

A first study of unemployment in France (Gaubert & Cottrell, 1999) has produced a broad classification of the unemployed. This classification is built using quantitative variables characterizing the unemployed (duration of unemployment, number of months of unemployment spent doing occasional work, job seniority, number of job offers received per month of unemployment).

Practicing an occasional work (AR) while being searching a durable job has been found to be the major discriminating factor in the studied population, while producing three broad classes (a fourth one is constituted by those looking for their first job, obviously different from the rest of the population):

- those looking for a durable job and not having had any occasional work during their search;
- those who, while searching a durable job, are having a part time activity which constitute a significant proportion of the duration of the spell of unemployment ; their situation seems to be closer to that of people having a precarious part-time occupation than to the unemployed according to the usual definition;
- those having had very few occasional jobs.

---

[1] We thanks ANPE and UNEDIC for a financial support for this study and the permission to use the data.

Simultaneously three types of exits from unemployment, using the usual definition, have to be distinguished: durable job, part time job with a limited duration and withdrawal from the labor market.

One of the characteristics of those having a precarious part time job is that they are looking for a durable job after the completion of a fixed-term contract. For them the duration of a spell of unemployment is significantly shorter than on the average. An interpretation may be that those employed under short-term employment contracts are having successive spells of full-time and part-time employment, with the same type of activity, possibly in the same firm. This could illustrate, on the individual's side, some effects of the firm's search for more flexibility.

This result would be an interesting one interpreted in terms of labor market segmentation, showing the major difficulty to find again a job in the primary segment after having unemployed.

Starting with these previous results, one purpose of this paper is to study the possible trajectories of the unemployed, emphasizing the successive observed transitions.

The data concerns unemployed people with two or more spells of unemployment, including information about the allowance received (periods and corresponding amounts) and the former wage in the last job occupied before the spell of unemployment.

One possible direction is to study the duration of the transition, the types of successive transitions followed, the role of the way unemployment benefit is designed in France (with a progressive reduction of the amount over time) on the specific path followed, in order to evaluate the main factors which determine recurring unemployment.

This new work is then a study of the segmentation of the French labor market in a way not frequently found in recent literature

For instance, some analyses may be indicated : in terms of career differencing (Favereau *et al.*, 1991 ; Theodossiou, 1995), returns on human capital (Balsan *et al*., 1994), or studies oriented towards firm organization, efficiency wage (Oi, 1990 ; Albrecht & Vroman, 1992), or discriminating practices (Dickens & Lang,, 1985 ; Boston, 1990).

## 2. DATA AND VARIABLES.

The initial data is the complete register of the unemployed held by the ANPE; information on unemployment benefits and compensations in latest job is added from the data collected by UNEDIC. The studied period goes from July 1993 to August 1996. The population is constituted of all the unemployed who were looking for a job at the beginning of this period[2], or who became unemployed later (but before the end of August 1996); at the end of the period they are either unemployed or their status has changed in some way.

For the present study we use a 1% sample of the unemployed registered in the administrative region of Ile-de-France (Paris and suburbs) having two or more spells of unemployment (590 000 individuals on a population of more than 2 millions 167 000), those named recurring unemployed. The resulting data is constituted by 19 246 individuals and 44 405 spells of unemployment.

This sample has been matched with the unemployment benefit register to retrieve the corresponding information.

Three types of information are available:

 - individual characteristics: age, marital status, number of children, nationality, level of education, seniority in the latest job, professional skill level;
 - characteristics of the unemployment: reason of job loss, how the unemployment status has ended, number of spells and duration of each one, number of job offers received, characteristics of the casual job if they has been some (number of hours per month and monthly wage)

---

[2] This explains that a duration of unemployment greater than 37 months may be found.

- characteristics of the unemployment allowance: average amount on a daily basis, compensation received in the latest job.

An *occasional job* (AR) is defined as any activity occupied by an individual who receives a compensation while he is registered as looking for a durable job. In France, this does not put an end to being registered as unemployed if the activity is very occasional or the number of hours per month remains small. Since 1995, the ANPE distinguishes among the registered unemployed those working for more than 78 hours per month. They constitute the so-called 6th category, excluded from some official statistics of unemployment, even though the person remain registered as looking for a job. In terms of unemployment benefits, if the characteristics of the latest job comply with the rules of unemployment benefits, they may perceive a compensation while its amount is not greater than 70 % of the latest job compensation. This occasional activity will be included to compute the next unemployment benefit.

The rules of computation of the *unemployment allowance*, on a daily basis, are quite complex. The rate is decreasing with thresholds defined differently according to the age and the professional history of the individual. For instance a young unemployed is submitted to a reduction rate of his benefits of 17 % each period of 122 days. In order to obtain a better homogeneity with the latest job compensation and the one obtained in an occasional activity, we have calculated an average allowance has been computed over the spell duration. This simplification removes the influence of the diminishing rate. We intend to investigate more thoroughly that point in a another study will.

The first step is a *simplification of the raw data* by two means. First, using analysis of variances and simple regressions, some variables have been removed as not significant (marital status, number of children, nationality, for instance). Then various qualitative variables have been coded with a number of modalities severely reduced. For instance the reasons of job loss are grouped in 4 classes instead of 11 in the administrative table, and 4 classes of exits from unemployment instead of 63.

The following tables summarize the **variables used**, with their name, definition and the codification for qualitative variables.

| Quantitative Variables | |
|---|---|
| AGE | Age of the individual |
| CMDUR | Cumulated duration of unemployment since initial registration in months |
| CPPAR | Proportion of the duration of the occasional work within the period of unemployment |
| DUR | Length of the latest period of unemployment in days |
| EXPER | Job seniority in years |
| INDUR | Cumulated duration of unemployment allowance in days |
| MGAIN | Monthly wage received in an occasional work |
| MXMHEUR | Amount of hours per month worked in an occasional job |
| NCHOM | Number of periods of unemployment |
| SRREVAL | Daily wage perceived in the latest job |
| TINDMOY | Average amount of daily benefits computed on the cumulated duration of unemployment |

**Table 1**

| Qualitative variables | |
|---|---|
| AGEC | Sub-categories of age: <25, 25-35, 35-45, 45-55, >55 |
| CTINDMOY | Daily benefits: <60 F, 60-100, 100-150, >150 |
| DIPL3 | Level of education: > bac (post secondary school level), bac level (secondary school completed), <bac (secondary school not completed) |
| DURC | Cumulated duration of unemployment: <12 months, 12-24, >24 |
| HAR | Monthly hours of occasional work (AR) :0, 0-39, 39-78, 78-117, >117 |
| PPARC | Proportion of cumulated duration of unemployment doing AR: 0, 0-0.1, 0.1-0.3, >0.3 |
| RMOTIFA | Types of exits from unemployment (4 categories detailed below) |
| RMOTIFI | Reasons for unemployment (4 categories detailed below) |

**Table 2**

The different types of *exits from unemployment* have been grouped in 4 categories:
1. job found by the individual himself or with the help of ANPE services;
2. enter in training program;
3. withdrawal from the labor market (illness, retirement, military service);
4. administrative cancellation (discouraged worker)[3].

Similarly, the causes of *registration at the ANPE* have been coded in 4 modalities:
1. lay off
2. end of fixed-term contracts
3. voluntary quit
4. first job search.

## 3. TRANSITIONS BETWEEN TWO STATUS.

The two next tables show the figures obtained, when the types of exits from unemployment are crossed with the causes of the former registration and, conversely, the causes of further registration crossed with the types of exits from current unemployment. Each cell shows the figure obtained by dividing the frequency of a particular transition by the total number of transitions (first number) and the one obtained by dividing this frequency by the row sum (second number) which is the total number of a type of registration. For instance the transition coded (1,1), from lay-off to a new job, represents 20.73 % of all the transitions and 59.08 % of those coming from lay-offs (which is itself 35.09 % of the different transitions).

| RMOTIFA RMOTIFI | 1 | 2 | 3 | 4 | Total |
|---|---|---|---|---|---|
| 1 | **20.73 %** **59.08 %** | 1.21 % 3.45 % | 4.81 % 13.71 % | **8.33 %** **23.75 %** | 35.09 % |
| 2 | **29.33 %** **62.67 %** | 1.57 % 3.35 % | 5.95 % 12.72 % | **9.95 %** **21.26 %** | 46.81 % |
| 3 | 4.94 % 70.29 % | 0.18 % 2.50 % | 0.86 % 12.31 % | 1.05 % 14.90 % | 7.03 % |
| 4 | **8.36 %** **75.51 %** | 0.41 % 3.72 % | 0.98 % 8.87 % | 1.32 % 11.90 % | 11.07 |
| Total | 63.37 % | 3.37 % | 12.61 % | 20.65 % | 100 % |

**Table 3 : Registration towards exit**

---

[3] Among the types defined by the ANPE corresponding to the last main category, class 90, involving a major difficulty. This class corresponds to people not complying the formal rules to remain registered; according to the ANPE services, most of these people have found a job by their own means and neglect to send an advice to the ANPE. So this class has finally been included in the first category.

It is noticeable that 5 types of transition are very significant. They are presented in bold characters, and coded (1,1), (1,4), (2,1), (2,4), (4,1). Respectively these are the following transitions
- lay-off ==> job
- lay-off ==> discouraged worker
- end of fixed-term contract (CDD) ==> job
- end of fixed-term contract ==> discouraged worker
- first job search ==> job.

The main weakness of the results presented in table 3 is that there is nothing about the kind of job found for the four cells of the first column. This is due to the fact that the information is not available in ANPE register, so it is impossible to separate durable jobs from fixed-term contracts (CDD).

The table presenting the further transitions gives a more precise information. People we refer to here are those who become unemployed again after an exit from the first spell of unemployment.

| RMOTIFI RMOTIFA | 1 | 2 | 3 | 4 | Total |
|---|---|---|---|---|---|
| 1 | **27.08 %** **34.62 %** | **41.20 %** **52.66 %** | 4.58 % 5.86 % | 5.37 % 6.86 % | 78.23 % |
| 2 | 0.86 % 48.22 % | 0.68 % 38.22 % | 0.06 % 3.56 % | 0.18 % 10.00 % | 1.79 % |
| 3 | **7.83 %** **56.11 %** | **4.86 %** **34.79 %** | 0.55 % 3.93 % | 0.72 % 5.18 % | 13.96 % |
| 4 | 3.05 % 50.76 % | 2.21 % 36.75 % | 0.27 % 4.49 % | 0.48 % 8.00 % | 6.01 |
| Total | 38.83 % | 48.95 % | 5.47 % | 6.75 % | 100 % |

**Table 4 : Transition matrix: latest exit to new registration**

Four types of transition are very significant; they are in bold characters and are coded (1,1), (1,2), (3,1), (3,2). Respectively they constitute the following transitions:
- job ==> lay-off
- job ==> end of fixed-term contracts of employment (end CDD)
- exit from the labor market ==> lay off
- exit from the labor market ==> end of fixed-term contracts (end CDD)

It is important to notice that the major transition to unemployment results from fixed-term contracts coming to an end. This means that a path constituted of consecutive periods of employment and unemployment is going to be the reference situation. The figures would be even greater if computed with data more recent than 1996.

With these transition tables, two new quantitative variables have been defined. RSS11 measures the proportion of transitions from job to lay-off among the total number of transitions 'off then back into unemployment'. RSS12 is obtained by dividing the number of transitions 'find a job then end of fixed-term contract' by the total number of transitions exit-registration. These variables are not used to perform the classification but they add some information to help the interpretation of the classes obtained.

## 4. CLASSIFICATION OBTAINED USING QUANTITATIVE VARIABLES.

A Kohonen classification have been constructed using the quantitative variables described in section 2. A Kohonen classification (Kohonen, 1984, 1993, 1995, Cottrell, Fort & Pagès, 1998) is a neural technique (Ripley, 1996) similar to the moving centers method (Lebart, Morineau & Piron, 1995) but using a notion of neighborhood between the classes. This implies a very useful property: neighbors are associated to the same class (as with any classification method) or in close classes. This produces a quantization of the data space which preserves the topology. This property allows very illustrative graphical representations (Kaski, 1997, Cottrell & Rousset, 1997, Cottrell, 1997, Cottrell et al., 1999) where the closeness between individuals is evaluated in a global way. The map obtained may be interpreted as a non linear projection of the data. The study of centroids (or code-vectors) of each class conducts to a clear interpretation of the classification. The successive steps are to begin with a detailed classification (100 classes in this study), then these classes are grouped using a clustering classification, to produce a small number of broad classes constituted with neighbor classes. The description of the broad classes defines a typology of the individuals.

In order to obtain the Kohonen classification, 10 quantitative variables are used (those listed in table 1 with the exception of the salary associated to the last job (SRREVAL) which is (by definition) not available for young people looking for a first job). From 100 classes at the beginning, 5 broad classes are eventually constituted for an easier description.

Two code-vectors are presented in Figure 1 and Figure 2 (selected from two very distant Kohonen classes). The line joining the values of the 10 variables represents the profile of this class. All the variables used are standardized over the whole observations in order to eliminate a level effect due to different magnitudes of the variables.

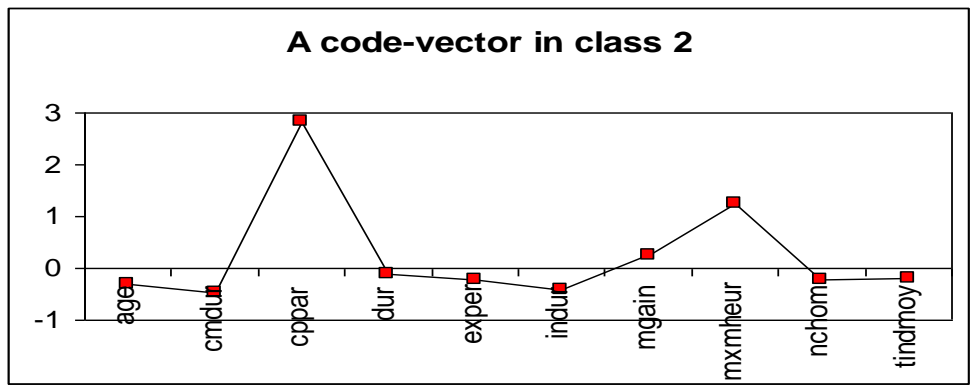

**Figure 1 : Code vector in class 2**

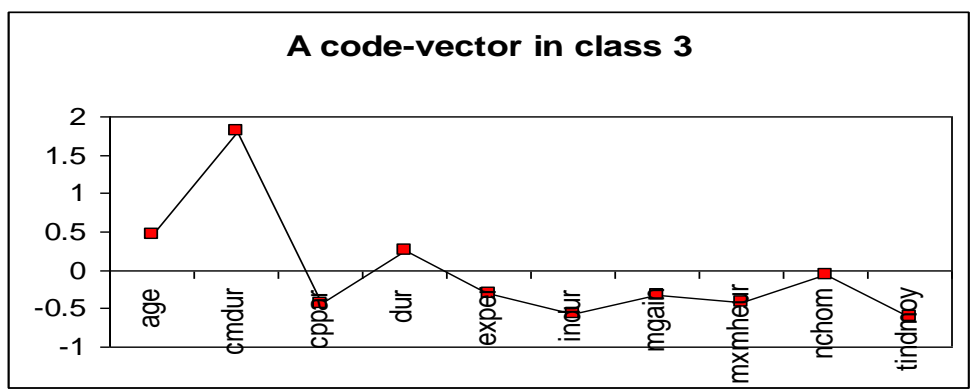

**Figure 2 : Code vector in class 3**

It is clear that the broad classes 2 and 3 are very different. Figure 3 shows the code vectors of the 100 classes and the grouping producing the broad classes, represented by different colors. The progressive transformation of the profiles may be noted.

Each of the broad classes is very different from the others by the number of individuals attached to each one.

| Classe 1 : 7908 | Classe 2 : 3793 | Classe 3 : 4519 | Classe 4 : 877 | Classe 5 : 2149 | Total : 19246 |
|---|---|---|---|---|---|

**Table 5 : Classes Frequencies**

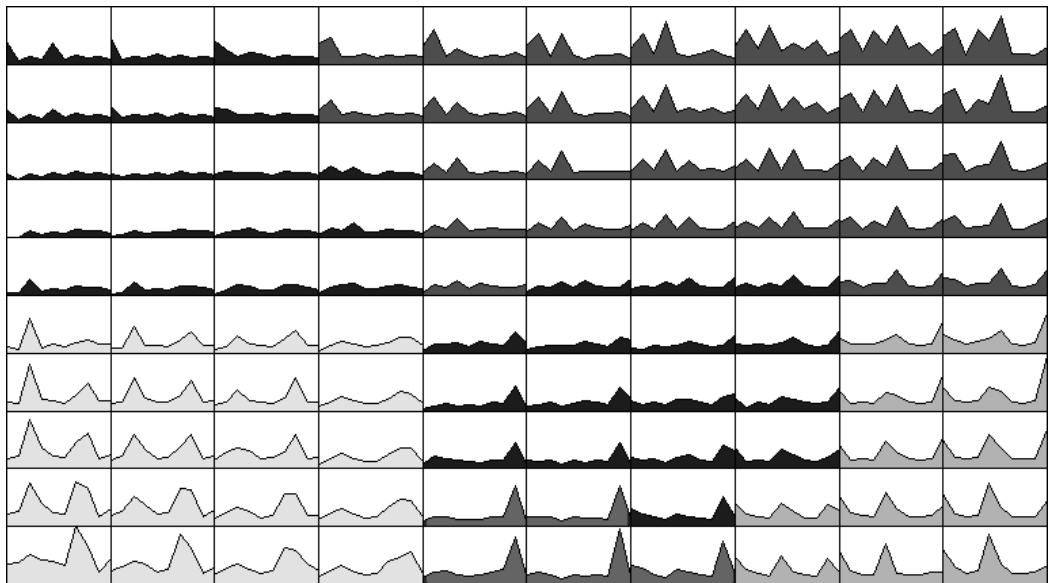

**Figure 3 : Kohonen Classification, classes and broad classes**
**(1 : Upper left corner and center, 2 : Down left corner, 3 : Upper right corner,**
**4 : Bottom center area, 5 : Bottom right corner)**

Figure 4 shows the inter-classes distances (Cottrell & de Bodt, 1996).

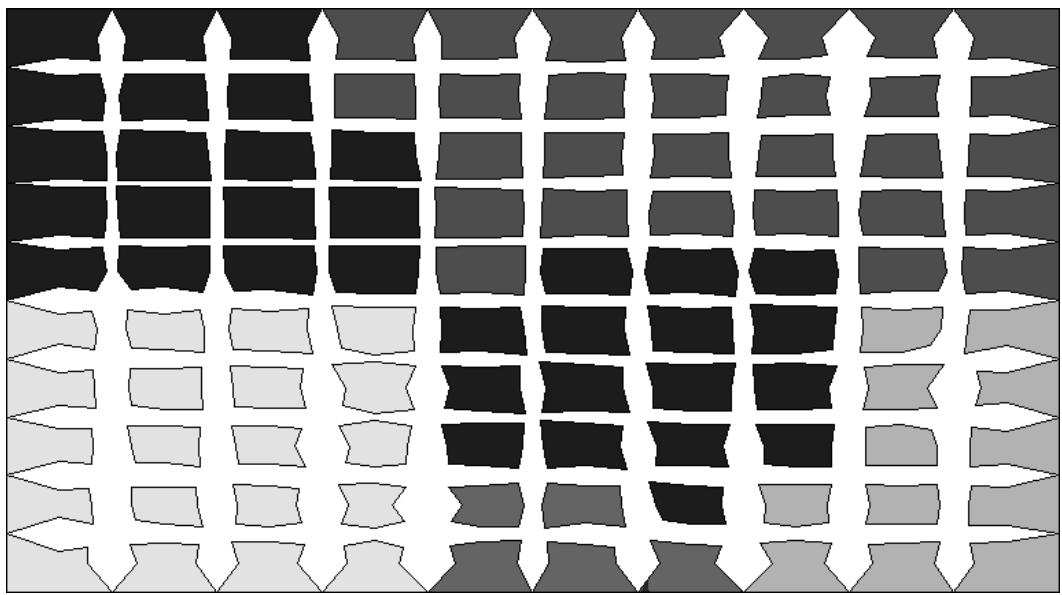

**Figure 4 : Inter-classes Distances and broad-classes**

In order to characterize the five broad classes, the means and standard deviations are examined (see table 6).

| | Classe 1 | Classe 2 | Classe 3 | Classe 4 | Classe 5 | Population |
|---|---|---|---|---|---|---|
| AGE Age of the individual | 28.78 (8.01) | 29.82 (8.06) | 34.80 (9.25) | 29.42 (8.16) | 44.78 (8.03) | 32.22 (9.75) |
| CMDUR Cumulated duration of unemployment since initial registration in months | 12.94 (9.80) | 24.27 (19.84) | 39.01 (20.00) | 22.23 (13.24) | 23.61 (17.63) | 22.91 (18.86) |
| CPPAR Proportion of the duration of the occasional work (AR) within the period of unemployment | 0.03 (0.06) | 0.31 (0.20) | 0.04 (0.08) | 0.05 (0.09) | 0.05 (0.09) | 0.09 (0.15) |
| DUR Length of the latest period of unemployment in days | 116.47 (96.95) | 281.94 (200.50) | 424.80 (253.95) | 112.97 (114.94) | 192.00 (155.37) | 229.75 (213.87) |
| EXPER Job seniority in years | 2.51 (3.36) | 3.81 (5.31) | 4.15 (5.25) | 2.88 (4.05) | 16.79 (8.82) | 4.76 (6.70) |
| INDUR Cumulated duration of unemployment allowance in days | 82.42 (131.35) | 169.03 (202.03) | 438.93 (414.16) | 121.49 (182.21) | 332.64 (245.01) | 212.92 (291.48) |
| MGAIN Monthly wage received in an occasional work (AR) | 29.07 (254.73) | 4612.14 (4386.92) | 334.47 (1124.81) | 126.66 (774.76) | 204.58 (944.29) | 1028.05 (2722.15) |
| MXMHEUR Amount of hours per month worked in an occasional job (AR) | 2.01 (9.53) | 127.88 (51.14) | 19.04 (44.06) | 11.92 (32.25) | 13.32 (36.32) | 32.53 (59.00) |
| NCHOM Number of periods of unemployment | 2.22 (0.42) | 2.26 (0.52) | 2.13 (0.34) | 4.36 (0.69) | 2.25 (0.49) | 2.31 (0.63) |
| TINDMOY Average amount of daily benefits computed on the cumulated duration of unemployment | 44.61 (59.80) | 79.99 (80.45) | 79.22 (63.52) | 44.00 (58.51) | 202.64 (155.71) | 77.32 (93.80) |
| SRREVAL Daily wage perceived in the latest job | 215.52 (124.31) | 263.19 (215.08) | 223.60 (191.18) | 230.71 (191.50) | 439.40 (325.71) | 265.70 (224.38) |
| RSSE11 Proportion of transitions out of-back into unemployment from job to lay off (1,1) | 0.27 (0.42) | 0.27 (0.42) | 0.24 (0.42) | 0.29 (0.33) | 0.32 (0.44) | 0.27 (0.42) |
| RSSE12 Proportion of transitions exit-registration from 'find a job' to end of contract (1,2) | 0.39 (0.47) | 0.45 (0.47) | 0.37 (0.47) | 0.46 (0.36) | 0.42 (0.47) | 0.40 (0.47) |

**Table 6 : Quantitatives Variables according to the 5 classes**

A more complete description of the classification may be obtained by using eight qualitative variables. Some of them are obtained by discretization of quantitative variables used to produce the classification (AGEC, CTINDMOY, DURC, HAR, PPARC). Others are exogenous (DIPL3, RMOTIFA, RMOTIFI).

## 5. TYPOLOGY OF CLASSES.

Figure 5 shows a representation of 8 qualitative variables for each of the broad classes and for the whole population.

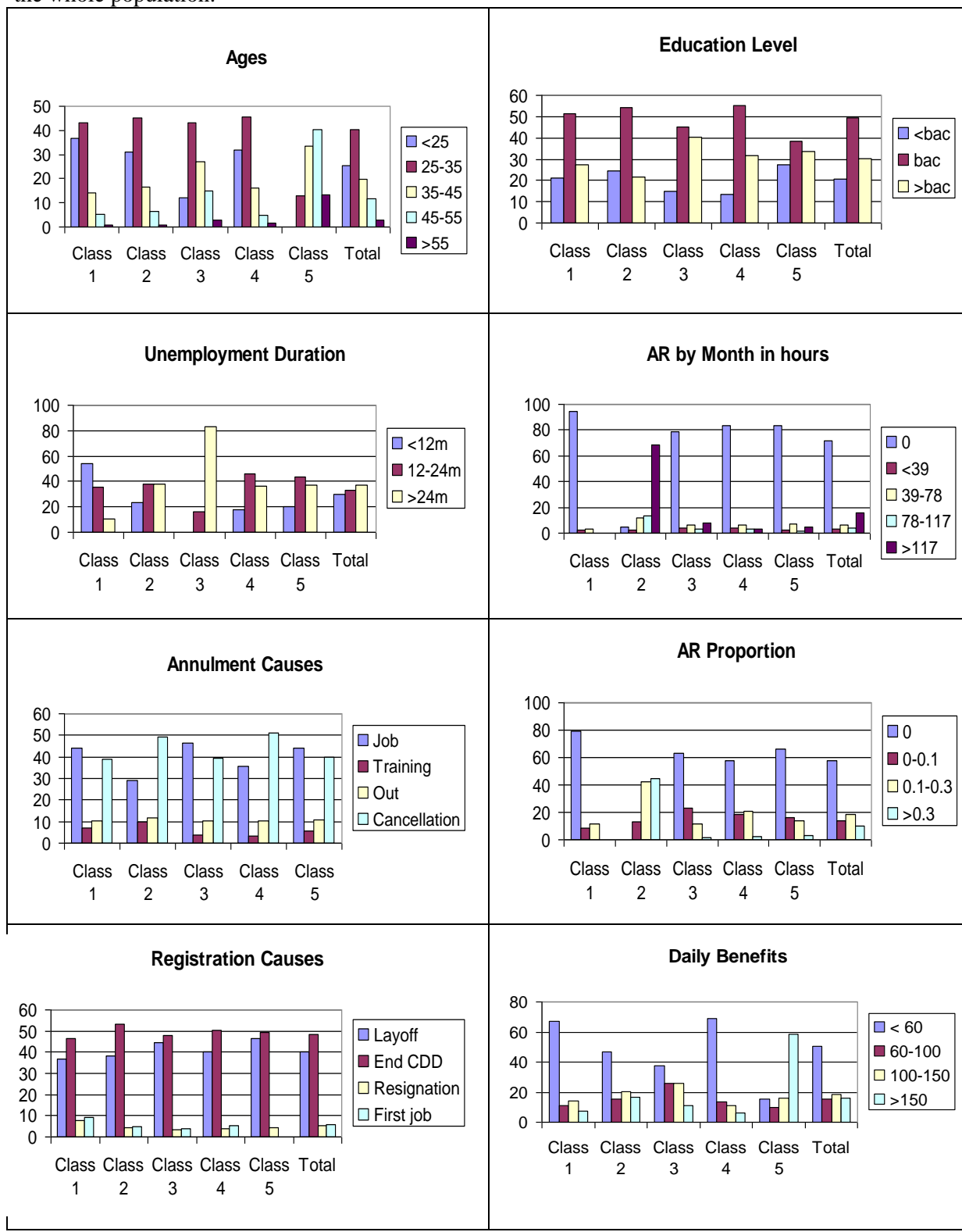

**Figure 5 : Qualitative Variables Distribution**

Class 1 is made up of young people, having a very short seniority, a duration of unemployment lower than the whole population mean, an average amount of unemployment benefits close to 0. Most of them are still looking for their first job.

Class 2 comprises people highly involved in occasional work. The duration of their unemployment and the benefits obtained are slightly above the average. The main cause of this unemployment is the end of fixed-term contracts.

Class 3 is characterized by the very long duration of unemployment. People are older than the average with benefits from unemployment slightly above the average. Most of them are not exerting any occasional work, except a sub-group which may be seen on figure 6.

Class 4 is constituted with young people, with no occasional work, and a short seniority in past employment. They are different in some points from those defined as class 1: the cumulated duration is longer, the greater recurrence of unemployment is greater, as the relative importance of ends of fixed-term contracts. Their situation seems a good illustration of a of a typical trajectory with successive periods of fixed-term contracts and unemployment.

Class 5 is made up of older people, with a lengthy experience and a corresponding higher averaged amount of benefits, rarely exerting an occasional work. The duration of their unemployment is slightly above the average but clearly lower than the one of class 3.

Thus it appears that there is no direct link between unemployment duration and the involvement in an occasional work or the amount of benefits perceived while unemployed. These results are clearly inconsistent with some usual conclusions applied to the unemployed as a whole. One of thoseis the presentation of various types of precarious jobs as a step leading to more durable jobs. Another one is the traditional positive link between the amount of benefits and the duration of unemployment. It appears that these conclusions are not verified for several subgroups of the unemployed.

More precisely it seems necessary to distinguish two types of occasional jobs. One with both an important amount of hours per month and a large proportion of the spell of unemployment occupied by this activity. This is a pattern very similar with that of precarious jobs, combining fixed-term employment contract and part time work. This type is represented here by class 2; the duration of unemployment is somewhat shorter than observed for people having some experience of the labor market[4]. One other type is constituted with people grouped in the classes 4 or 5, rarely exerting this kind of activity and with an amount of monthly hours frankly lower than the one of the first type. Their duration of unemployment is slightly longer than the population average.

On the subject of the effect of unemployment benefits on unemployment duration it is possible that the reduction resulting from the use of calculated average amounts perceived results in an important attenuation of this effect. This may be of some importance for the first months of a spell of unemployment, when the diminishing rate used to define the amount of benefit implies for the unemployed a real reduction of the perceived amount.

Another general observation is that the obtained classes are not differentiated when considering the transition matrix 'exit from unemployment then back to unemployment': the main transition is clearly the end of a fixed-term contract found after a spell of unemployment.

This analysis may be improved with a study of the distribution of the population observed using qualitative variables : the distribution of each one within the Kohonen classes and its transformation when describing the whole grid. This transformation may induce specific groupings of classes, possibly with an easier interpretation. An example is produced by Figure 6 representing the different types of exercise of an occasional job. Another one is Figure 7 representing the different levels of daily benefits.

---

[4] Thus excluding those looking for their first job characterized with a significant shorter duration of unemployment.

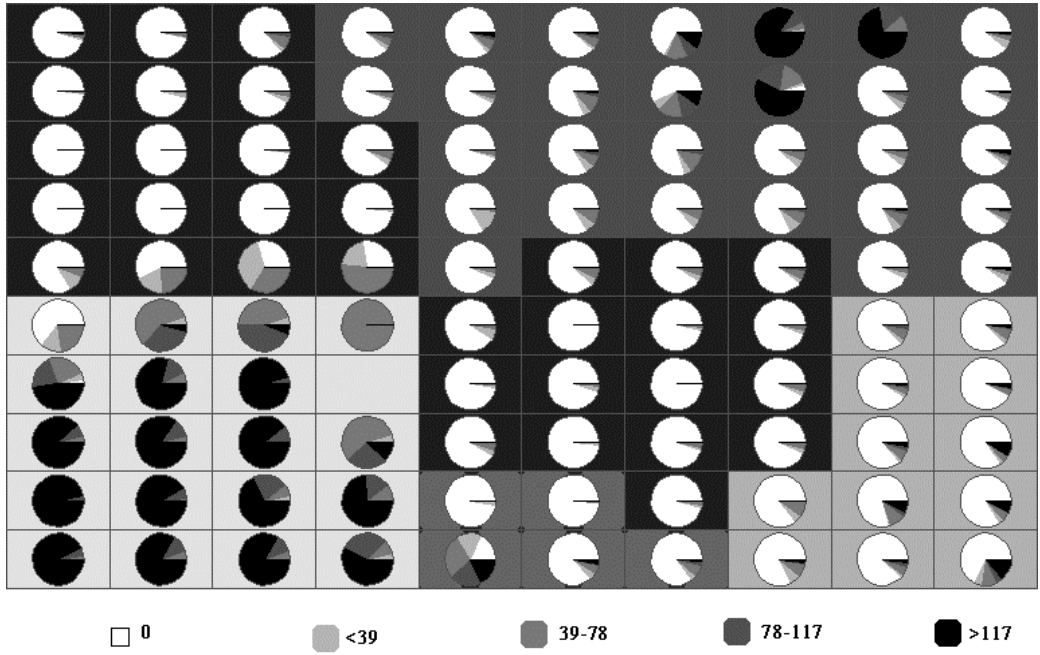

**Figure 6 : Distribution of Monthly Hours of Occasional Work (AR) in the broad classes**

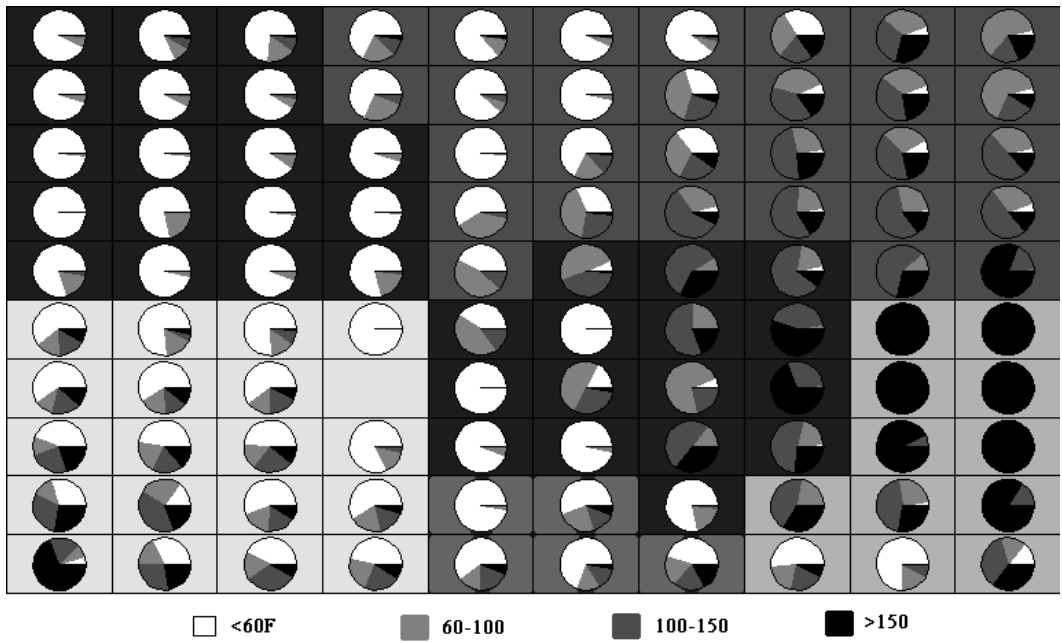

**Figure 7 : Distribution of the Daily Benefits in the broad classes**

## 6. ANALYSIS OF MULTIPLE CORRESPONDENCES.

A complete description of the characteristics of the unemployed and their trajectories is conducted using the qualitative variables presented earlier plus a new one defined as the broad class in which each individual is placed. A classical factor analysis is carried out (Lebart, Morineau & Piron, 1995), the multiple correspondence analysis. It used 8 qualitative variables, representing 32 modalities.

The following figure is the projection of the set of modalities obtained using the two first dimensions.

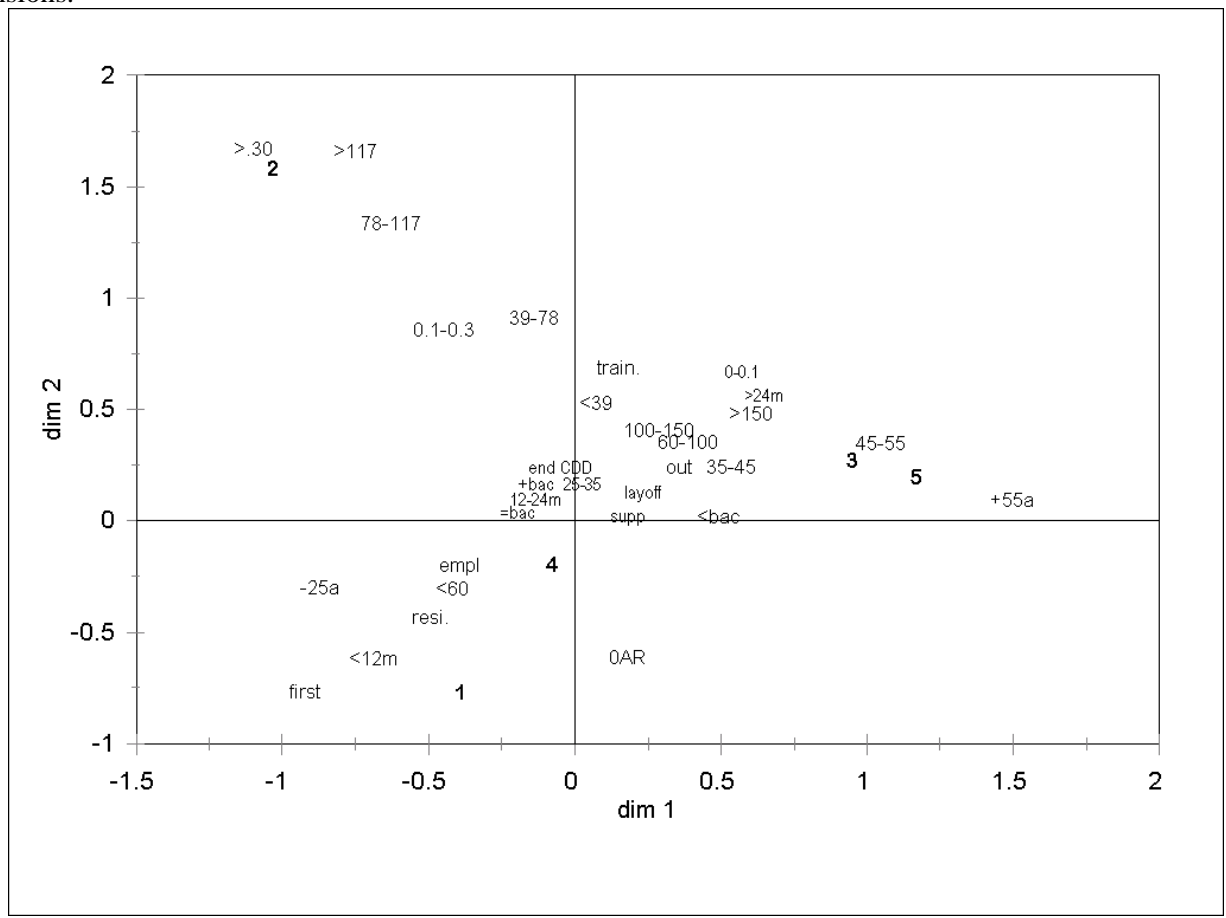

**Figure 8 : The two first axes.**

It shows some results already found in the previous steps. The first axis is clearly defined with the amount of hours corresponding to an occasional job.

The second one shows a contrast between young people with a high educational level and older people with an educational level below from the average. This axis is very close to the one defined by the duration of unemployment, causes of unemployment and types of exits from unemployment. Young people are associated with shorter duration of unemployment, search for a first job and exit from unemployment with a job. Causes of unemployment are placed along this axis with the following order, first job, quit, end of fixed-term contract, lay-off. The types of exit from unemployment are similarly projected along the same axis showing the exit by a job sharply separated from the three other ones, almost stacked at the opposite side of the axis.

Involvement in an occasional work and the whole set of other factors are orthogonal, as it has been observed earlier.

The projection on the axes 1 and 3 leads to the same conclusions. See Figure 9.

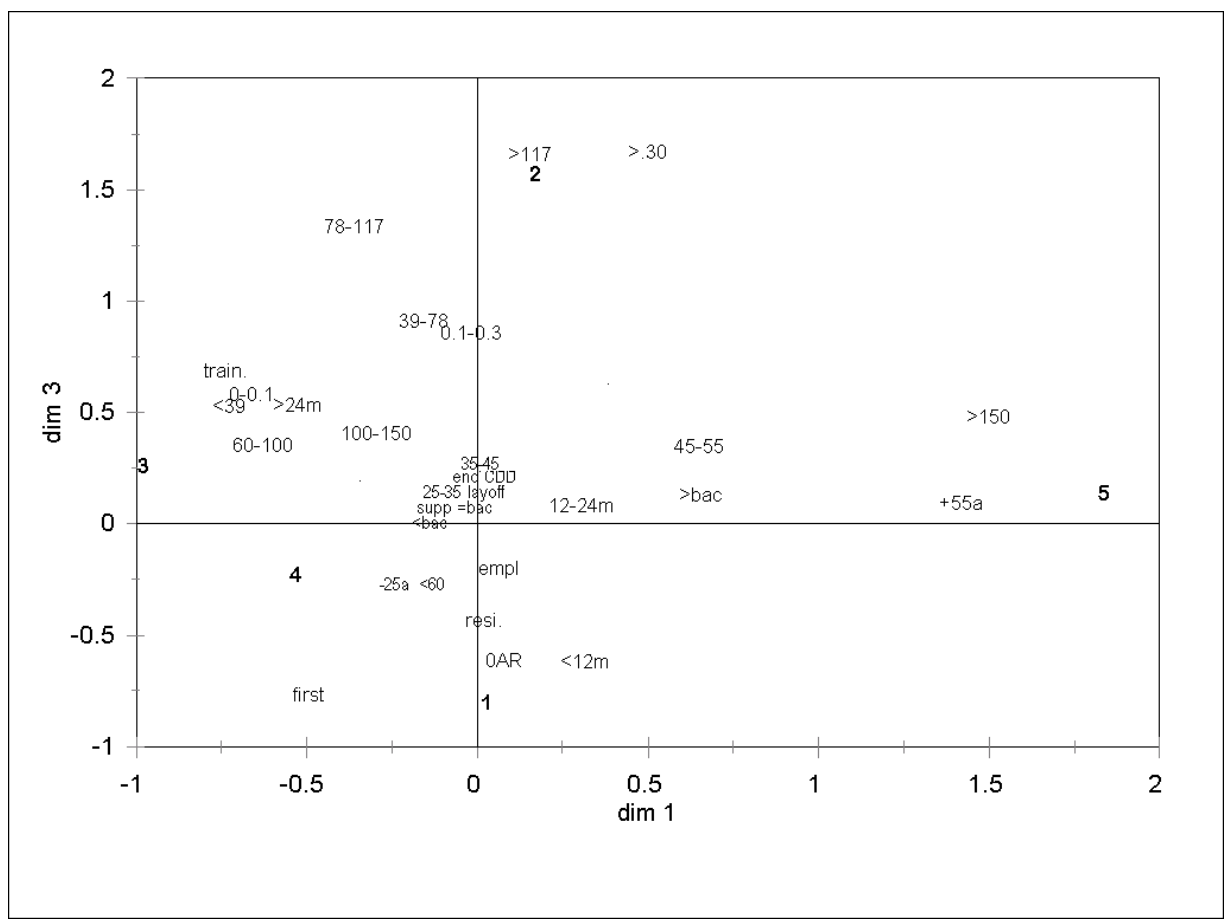

**Figure 9 : The axes 1 and 3**

## 7. CONCLUSION.

A typology of recurring unemployed may be defined. It is based on age, level of education, experience, involvement in an occasional work, amount of benefits from unemployment. But it seems that there is no link between the practice of an occasional work and the duration of unemployment. About the question of unemployment recurrence, the main fact is that the exit from unemployment is typically obtained with a fixed-term job, which in turn determines the next unemployment spell. The observed cumulated duration of unemployment is 23 months for the whole population studied and respectively, for the broad classes, 13, 24, 39, 22 and 23 months.

A further step is to analyze the impact of unemployment benefits and the diminishing rate applied in a more precise way. This means a study of the successive periods constituting a spell of unemployment and the possible impact of the decreasing amount of allowance on the possible exit from unemployment.

This study would be largely improved with a thorough analysis of the type of occupation attached to the occasional job and its possible correspondence to the one of the regular job. This is an important feature to understand the phenomenon of recurrence. The problem is that, at this time, this information is not available in the data furnished by the various institutions involved.

From a technical point of view, it could be very useful to apply another neural technique called KACM (Ibbou, 1998, Cottrell, et al., 1999), in order to obtain the modalities of qualitative variables on a unique grid, instead of the successive planes furnished by the traditional factor analysis.